\documentclass[10pt]{article}

\usepackage{latexsym,amsmath,graphics,graphicx}
\usepackage[psamsfonts]{amssymb}
\usepackage[latin9]{inputenc}
\usepackage[T1]{fontenc}
\usepackage[bookmarks,bookmarksnumbered,pdfstartview=FitBH]{hyperref}
\usepackage{color}
\definecolor{darkred}{rgb}{0.5,0,0}
  \definecolor{darkgreen}{rgb}{0,0.5,0}
  \definecolor{darkblue}{rgb}{0,0,0.5}

\hypersetup{
colorlinks=true, 
linkcolor= darkblue, 
citecolor= darkgreen, 
urlcolor= blue,
}
\newtheorem{theorem}{Theorem}[section]

\newtheorem{lemma}[theorem]{Lemma}

\newtheorem{proposition-def}[theorem]{Proposition-Definition}
\newtheorem{corollary}[theorem]{Corollary}
\newtheorem{definition}[theorem]{Definition}

\def\RR{\mathcal{R}}
\def\SS{\mathcal{S}}

\def\FF{\mathbf{F}}
\def\Nmath{\mathbb{N}}
\def\Zmath{\mathbb{Z}}
\def\Rmath{\mathbb{R}}

\def\Tmath{\mathbb{T}}

\def\HHH{\mathcal{H}}
\def\point{\mathop{\cdot}}
\def\WW{\mathcal{W}}

\def\OE{\raisebox{-0.5mm}{\ ${{\overset{\mathrm{OE}}{\sim}}}$\ }}
\def\Conj{\raisebox{-0.5mm}{\ ${{\overset{\mathrm{Conj}}{\sim}}}$\ }}

\def\SL{\mathrm{SL}}
\def\FFF{\mathcal{F}}

\newcommand{\acting}[1]{\overset{#1}{\curvearrowright}}

\def\freeprod{\mathop{*}}

\def\dom{\mathrm{dom}}
\def\rng{\mathrm{rng}}

\def\actionorig{\sigma}

\def\ll{\langle}
\def\rr{\rangle}

\def\Out{\mathrm{Out}}
\def\Aut{\mathrm{Aut}}
\def\EEE{\mathcal{E}_{1}}
\def\RRTE{\widetilde{\RR}_E}
\def\RRT{\widetilde{\RR}}
\def\RRTF{\widetilde{\RR}_F}
\def\RRE{\RR_{E}}

\begin{document}

\thispagestyle{empty}
\title{Free Product Actions with Relative Property (T) and Trivial Outer Automorphism Groups}
\author{Damien Gaboriau\thanks{CNRS} }
\date{November 18, 2009}
\maketitle

\begin{abstract}
{We show that every non-amenable free product of groups admits free ergodic probability measure preserving actions which have relative property (T) in the sense of S.~Popa \cite[Def. 4.1]{Pop06}. There are continuum many such actions up to orbit equivalence and von Neumann equivalence, and they may be chosen to be conjugate to any prescribed action when restricted to the free factors.
We exhibit also, for every non-amenable free product of groups, free ergodic probability measure preserving actions whose associated equivalence relation has  trivial outer automorphisms group.
This gives, in particular, the first examples of such actions for the free group on $2$ generators.}
\end{abstract}


\section{Introduction}

Several breakthroughs in von Neumann Algebras theory and Orbit Equivalence have been made possible, during the last years, by the introduction by S. Popa of the notion of \textbf{rigidity} for pairs $B\subset N$ of von Neumann algebras \cite{Pop06}. This property is inspired by the relative property (T) of Kazhdan and is satisfied by the inclusion $L^{\infty}(\Tmath^2)\subset L^{\infty}(\Tmath^2)\rtimes \SL (2, \Zmath)$ associated with the standard action of $\SL(2,\Zmath)$ on the $2$-torus. It has proved to be extremely useful, more generally, for inclusions $A\subset M(\RR)$ arising from standard countable probability measure preserving (p.m.p.) equivalence relations $\RR$, where $A$ is a Cartan subalgebra  in the generalized crossed product or group-measure-space construction von Neumann algebra \cite{MvN36}, \cite{FM77b}.

 In case the pair $A\subset M(\RR)$ is rigid, the relation is said to have the \textbf{relative property (T)}. When combined with antagonistic properties like Haagerup property \cite{Pop06} or (amalgamated) free product decomposition \cite{IPP05}, the rigid Cartan subalgebra were shown to be essentially unique, thus allowing orbit equivalence invariants (like the cost \cite{Gab00a}, $L^2$-Betti numbers $\beta_n(\RR)$ \cite{Gab02} or the fundamental group $\FFF(\RR)$ \cite[Cor. 5.7]{Gab02},\dots) to translate to von Neumann algebras invariants. 
This led to the solution of the long standing problem of finding von Neumann II$_1$ factors with trivial fundamental group \cite{Pop06}.

While the class of groups that admit a free p.m.p. ergodic relative property (T) action is closed under certain algebraic constructions (like direct products, commensurability \cite{Pop06},  or free product with an arbitrary group \cite[Cor. 7.15]{IPP05}), the building blocks were very few and relying on some arithmetic actions \cite{Pop06}, \cite{Val05}, \cite{Fer06} leaving open the general problem \cite[Prob. 5.10.2.]{Pop06}:
\emph{"Characterize the countable discrete groups $\Gamma_0$ that can act
with relative property (T)  on the probability space $(X,\mu)$, i.e., for which there exist free
ergodic measure preserving actions $\sigma$ on $(X, \mu)$ such
that $L^\infty(X, \mu) \subset L^\infty(X, \mu) \rtimes_\sigma
\Gamma_0$ is a rigid embedding."}

The purpose of this paper is first to prove that the class of groups 
that admit such a free p.m.p. ergodic relative property (T) action
contains all non-amenable free products of groups. Moreover, we show that they have continuously many different such actions (Th.~\ref{th: continuum many non vN act prescribed on factors})
and we remove any arithmetic assumption on the individual actions of the building pieces. In fact, given the state of art, an arithmetic flavor remains hidden in the way the individual actions are mutually arranged.

\bigskip

Outer automorphisms groups of standard equivalence relations are usually hard to calculate and there are only few special families of group actions for which one knows that $\Out(\RR)=\{1\}$. The first examples are due to S.~Gefter \cite{Gef93}, \cite{Gef96} and more examples were produced by A.~Furman \cite{Fur05}. They all take advantage of the setting of higher rank lattices. Monod-Shalom \cite{MS06} produced an uncountable family of non orbit-equivalent actions with trivial outer automorphisms group, by left-right multiplication on the orthogonal groups of certain direct products of subgroups.  Using rigidity results from \cite{MS06}, Ioana-Peterson-Popa \cite{IPP05} gave the first shift-type examples. However, all these examples are concerned with very special kind of actions. And the free groups keep out of reach by these techniques. 

The first general result appeared recently in a paper by S.~Popa and S.~Vaes \cite{PV08} and concerns free products $\Lambda*\FF_{\infty}$ of any countable group $\Lambda$ with the free group on infinitely many generators $\FF_{\infty}$, with no condition at all on the $\Lambda$-action. This relies heavily on the existence of a relative property (T) action for $\FF_{\infty}$.
The second purpose of our paper is twofold: extend the result from $\FF_{\infty}$ to free groups on finitely many generators, and using the first part extend it to any non elementary free products $\Gamma_1*\Gamma_2$ of groups (i.e. $(\vert \Gamma_1\vert -1)(\vert \Gamma_2\vert -1)\geq 2$).
This leads, in particular, to the first examples of actions of the free groups $\FF_n$, $\infty > n >1$, with trivial outer automorphisms group.

\begin{center}
{------------oooOOOOOooo------------}
\end{center}

Before stating more precisely our results, let's recall some definitions.
We only consider measure preserving actions on the standard Borel space.
First, from \cite[Def. 4.1]{Pop06} the definition of a \textbf{rigid inclusion} (or of relative rigidity of a subalgebra).
\begin{definition} \label{def: rigidity}
Let $M$ be a factor of type II$_1$ with normalized trace
$\tau$ and let $A \subset M$ be a von Neumann subalgebra.
The inclusion $A \subset M$ is called \textbf{rigid} if the following
property holds: for every $\epsilon > 0$, there exists a finite subset
$\mathcal{J} \subset M$ and a $\delta > 0$ such that whenever $_M H_M$ is
a Hilbert $M$-$M$-bimodule admitting a unit vector $\xi$ with the properties
\begin{itemize}
\item $\| a \cdot \xi - \xi \cdot a \| < \delta$ for all $a \in \mathcal{J}$,
\item $| \langle a \cdot\xi,  \xi \rangle - \tau(a)| < \delta$ and
$|\langle \xi \cdot a , \xi \rangle - \tau(a) | < \delta$ for all
$a$ in the unit ball of $M$,
\end{itemize}
then, there exists a vector $\xi_0 \in H$ satisfying $\|\xi - \xi_0\| < \epsilon$
and $a \cdot \xi_0 = \xi_0 \cdot a$ for all $a \in A$.
\end{definition}

\begin{definition}\cite[Def. 5.10.1]{Pop06}
A free p.m.p. ergodic action $\Gamma\acting{\sigma} (X,\mu)$ (respectively a countable standard p.m.p. equivalence relation $\RR$ on $(X,\mu)$) is said to have the \textbf{relative property (T)} if the inclusion of the Cartan subalgebra $L^{\infty}(X,\mu)$ is rigid in the (generalized) crossed-product  $L^{\infty}(X,\mu)\subset L^{\infty}(X,\mu)\rtimes_{\sigma} \Gamma$ (resp. $L^{\infty}(X,\mu)\subset M(\RR)$).
\end{definition}
One could say that the equivalence relation has ``the property (T) relative to the space $(X,\mu)$''.
Observe that A.~Ioana \cite[Th. 4.3]{Ioa07} exhibited for every non-amenable group, some actions $\sigma$ satisfying a weak form of relative rigidity, namely for which there exists a diffuse $Q\subset L^{\infty}(X,\mu)\subset M(\RR_{\sigma})$ such that $Q$ is relatively rigid in $M(\RR_{\sigma})$ and has relative commutant contained in $L^{\infty}(X,\mu)$. 

Recall the following weaker and weaker notions of equivalence for p.m.p. actions or standard equivalence relations $\RR, \SS$ on $(X,\mu)$:\\
Two actions $\Gamma\acting{\alpha}(X,\mu)$ and $\Lambda\acting{\beta}(X,\mu)$ are \textbf{Conjugate} \begin{equation}
\Gamma\acting{\alpha}(X,\mu)\Conj\Lambda\acting{\beta}(X,\mu)
\end{equation}
 if there is a group isomorphism $h : \Gamma\to \Lambda$ and a p.m.p. isomorphism of the space $f:X\to X$ that conjugate the actions $\forall \gamma\in \Gamma, \text{ (almost every) } x\in X$: 
$f(\alpha(\gamma)(x))=\beta(h(\gamma))(f(x))$.

Two p.m.p. standard equivalence relations $\RR, \SS$  are \textbf{Orbit Equivalent} 
\begin{equation}
\RR\OE\SS
\end{equation} 
if there is a p.m.p. isomorphism of the space that sends classes to classes, 
 or equivalently \cite{FM77b} if the associated pairs are isomorphic
\begin{equation}
\bigl(L^{\infty}(X,\mu)\subset M(\RR)\bigr) \ \ \simeq \ \ \bigl(L^{\infty}(X,\mu)\subset M(\SS)\bigr).
\end{equation}
This makes the relative property (T) an orbit equivalence invariant.

The standard equivalence relations 
 are \textbf{von Neumann Equivalent} if solely the generalized crossed products are isomorphic
\begin{equation}
M(\RR) \ \ \simeq \ \ M(\SS).
\end{equation} 

They are \textbf{von Neumann Stably Equivalent} if (they are ergodic and) the generalized crossed product factors are stably isomorphic for some amplification $r\in \Rmath^{*}_{+}$
\begin{equation}
M(\RR) \ \ \simeq \ \ M(\SS)^{r}.
\end{equation} 

We obtained in \cite{GP05} that the non-cyclic free groups admit continuously many relative property (T) orbit inequivalent, and even von Neumann stably inequivalent, free ergodic actions. 
We show that free products admit relative property (T) free actions and that we have a full freedom for the conjugacy classes of the
restrictions of the action to the free factors. In what follows,  a free product decomposition $\Gamma=\freeprod_{i\in I} \Gamma_i$ is called \textbf{non elementary} if $(\vert I\vert -1)\prod_{i\in I} (\vert \Gamma_i\vert -1)\geq 2$, i.e. there are at least $2$ free factors, none of them is $\sim\{1\}$ and if $\vert I\vert =2$, then one of the $\Gamma_i$ has at least $3$ elements.

\begin{theorem}\label{th: continuum many non vN act prescribed on factors}
Every non elementary free product  
$\Gamma=\freeprod\limits_{i\in I} \Gamma_i$
admits continuum many von Neumann stably inequivalent relative property (T) free ergodic p.m.p. actions,
whose restriction to each factor is conjugate with any (non necessarily ergodic) prescribed free action.
\end{theorem}
More precisely, let $\Gamma_i\acting{\actionorig_i}(X,\mu)$ be an at most countable collection of p.m.p. (non necessarily ergodic) free actions of countable groups $\Gamma_i$ on the standard probability space.
There exists continuum many von Neumann inequivalent free ergodic actions $(\alpha_t)_{t\in T}$ of the free product $\Gamma=*_{i\in I} \Gamma_i$ 
that have relative property (T), and such that for every $t\in T$ and $i\in I$, the restriction $\alpha_t\vert\Gamma_i$ of $\alpha_t$ to $\Gamma_i$ is conjugate with $\actionorig_i$
\begin{equation}
\bigl(\Gamma_i\acting{\alpha_t\vert \Gamma_i}(X,\mu)\bigr)\Conj\bigl(\Gamma_i\acting{\actionorig_i}(X,\mu)\bigr).
\end{equation}

\medskip

We introduced in \cite{Gab00a} in connection with cost, the notion of freely independent equivalence relations and of free decomposition of an equivalence relation (see \cite{Alv08} for a geometric approach):
\begin{equation}
\SS=\freeprod_{i\in I} \SS_i.
\end{equation}
The following, essentially due to A.~T{\"o}rnquist \cite{Tor06}, see also  {\cite[Prop. 7.3]{IPP05}}, states that countable standard equivalence relations may be put in general position:
\begin{theorem}[Tornquist]\label{th: S-i in general position}
Let $(\SS_i)_{i\in I}$ be a countable collection of standard p.m.p. equivalence relations on $(X,\mu)$, then there exists an equivalence relation $\SS$ on $(X,\mu)$ generated by a family of freely independent subrelations $\SS'_i$ such that $\SS'_i\OE \SS_i$, $\forall i\in I$. 
\end{theorem}
We obtain in fact the following more general form of Theorem~\ref{th: continuum many non vN act prescribed on factors} involving equivalence relations instead of free actions. Recall that a p.m.p. standard equivalence relation is called \textbf{aperiodic} if almost all its classes are infinite.
\begin{theorem}\label{th: continuum many non vN rel. prescribed on factors}
Let $(\SS_i)_{i\in I}$ be a countable collection (with $\vert I\vert \geq 2$) of p.m.p. standard countable aperiodic equivalence relations on the standard non atomic probability space $(X,\mu)$. 
Then there exists continuum many von Neumann stably inequivalent relative property (T) ergodic p.m.p. equivalence relations
on $(X,\mu)$ generated by a freely independent family of subrelations $\SS'_i$ such that for every $i\in I$, $\SS'_i\OE \SS_i$.\\
More precisely, there exists a strictly increasing continuum of ergodic equivalence relations $\SS_{t}$, $t\in (0,1]$ such that
\begin{enumerate}
\item for every $t\in (0,1]$, we have a free decomposition $\SS_{t}=\freeprod\limits_{i\in I} \SS_{i,t}$ such that $\SS_{i,t}\OE \SS_i$ for all $i\in I$
\item  for every $t\in (0,1]$,  $\lim\limits_{s\to t, s<t}\nearrow \SS_{s}=\SS_{t}$
\item $\SS_t$ has the relative property (T) 
\item the classes of stable isomorphism among the $M(\SS_{t})$ are at most countable,
in particular there is an uncountable set $T$ such that for $t\in T\subset (0,1]$  the $M(\SS_{t})$  are pairwise not stably isomorphic.
\end{enumerate}
\end{theorem}

\bigskip
S.~Popa and S.~Vaes obtained actions with trivial outer automorphism group for the particular family of free products $\FF_\infty * \Lambda$ (and $\Lambda$ may be trivial):
\begin{theorem}[{\cite[Th.~4.3]{PV08}}]
For every countable group $\Lambda$, there exists a continuum of pairwise von Neumann stably inequivalent free, ergodic, p.m.p. actions  $(\sigma_t)_{t \in T}$ of the free product $\FF_\infty * \Lambda$, with relative property (T), with $\Out(\RR_{\sigma_t})=\{1\}$, with trivial fundamental group, and whose restriction  to $\Lambda$ is conjugate with a prescribed free action.
\end{theorem}
Inspired by their techniques, we extend this kind of results to the non-cyclic free groups on finitely many generators, and more generally to any free product without any condition on the building blocks:
\begin{theorem}\label{th: F3-action with trivial out}
There exists continuum many, pairwise von Neumann stably inequivalent,
relative property (T), ergodic p.m.p. free actions $(\sigma_t)_{t\in T}$ of the free group $\FF_r$, $r=2, 3, \cdots$ such that
 $\Out(\RR_{\sigma_t})=\{1\}$.
 \end{theorem}
\begin{theorem}\label{cor: free prod gps out=1}
Every non elementary free product of countable groups $\freeprod_{i\in I} \Gamma_i$ admits continuum many von Neumann inequivalent 
free ergodic $(\sigma_t)_{t\in T}$ actions with relative property (T), with  $\Out(\RR_{\sigma_t})=\{1\}$ and whose restriction to each factor $\Gamma_i$ is conjugate with a prescribed free action.
\end{theorem}
This follows from the more general form involving equivalence relations.
\begin{theorem}\label{th: cont many non vN rel. Out=1 prescribed on factors}
Let $(\SS_i)_{i\in I}$, $\vert I\vert \geq 2$,  be a countable collection of p.m.p. standard countable aperidodic equivalence relations on the standard Borel space. 
Then there exists a continuum of pairwise von Neumann inequivalent,
ergodic equivalence relations $\SS_{t}$, with relative property (T) and $\Out(\SS_{t})=\{1\}$ such that:
\begin{equation*}
\SS_{t}=\freeprod\limits_{i\in I} \SS_{i,t} \textrm{\ \ \ \ \ for some \ \ \ \ \ } \SS_{i,t}\OE \SS_i \textrm{\ \ \ \ \ for all\ } i\in I \textrm{
 \ and for all\ } t\in T.
\end{equation*}
\end{theorem}

Since they satisfy property $\mathcal{FT}$, it follows from \cite{IPP05} that von Neumann stable equivalences between 
the von Neumann algebras $M(\SS_{t})$ entail stable orbit equivalences between the equivalence relations themselves, and with the same compression constant, in particular the fundamental group of $M(\SS_{t})$ coincides with the fundamental group of $\SS_t$. Thus, under mild conditions on $L^2$-Betti numbers (see \cite{Gab02}) or on the cost (see \cite{Gab00a}) we get a more precise information and produce plenty of new examples of factors with trivial fundamental group:
\begin{corollary}
If some $L^2$-Betti number $\beta_n(\SS_{t})$ for some $n\in \Nmath\setminus\{0\}$, or if the $\mathrm{cost}(\SS_t)-1$ does not belong to $\{0,\infty\}$  then the von Neumann inequivalent $\SS_t$ in Theorems~\ref{th: continuum many non vN act prescribed on factors} and \ref{th: cont many non vN rel. Out=1 prescribed on factors} are in fact von Neumann stably inequivalent $\SS_t$ and the associated von Neumann algebras 
have trivial fundamental group:
$\mathcal{F}(M(\SS_{t}))=\{1\}$.
This is also the case for the above $\RR_{\sigma_t}$ for free products of finitely many finitely generated groups.
\end{corollary}
Observe that the $\beta_n(\SS_{t})$ and $\mathrm{cost}(\SS_t)$ only depend on those of the free factors, namely:
$\beta_n(\SS_{t})= 
\sum_{i\in I}\beta_n(\SS_i)$ for $n\not=1$, and $\beta_1(\SS_{t})=\vert I\vert +\sum_{i\in I}\beta_1(\SS_i)$ and $\mathrm{cost}(\SS_t)=\sum_{i\in I} \mathrm{cost}(\SS_i)$.

\section{Proofs}

\textbf{Notation}

\begin{itemize}
\item $\RR\vert Y$ the restriction of $\RR$ to a Borel subset $Y\subset X$. 

\item $\FF_n$ the free group on $n$ generators

\item $\RR_{\beta}$ the equivalence relation associated with an action $\beta$.

\item $\langle \RR_1, \RR_2, \cdots, \RR_k\rangle$ the equivalence relation generated by $\RR_1, \RR_2, \cdots, \RR_k$.

\item $\langle \varphi_1, \varphi_2, \cdots, \varphi_n\rangle$ the equivalence relation generated by the family of partial isomorphisms $(\varphi_1, \varphi_2, \cdots, \varphi_n)$.

\item $\dom(\varphi)$ and $\rng(\varphi)$ denote the domain and range of the partial isomorphism $\varphi$

\item $\RR_1\freeprod \RR_2$ means that $\RR_1$, $\RR_2$ are freely independent, and $\RR=\RR_1\freeprod \RR_2$ that $\RR$ is generated by the freely independent subrelations $\RR_1$ and $\RR_2$

\item $\RR^p$ the $p$-amplification, i.e. the equivalence relation on $X\times\{1, 2, \cdots, p\}$ defined by $(x,j)\sim (y,k)$ iff $x\RR y$. 
Thus, identifying $X$ with $X\times \{1\}$, 
\begin{equation}\label{eq:amplification and free prod}
\RR^p= \RR*\langle \varphi_{2}\rangle*\cdots*\langle \varphi_{p}\rangle
\end{equation}
where  $\varphi_j:X\times\{1\}\to X\times \{j\}$ ($j=2, \cdots, p$) is the partial isomorphism defined by the identity on the first coordinate. 
\end{itemize}

If $\RR$ is ergodic, one defines  (up to isomorphism) the \textbf{compression} $\RR^{t}$ for $t\in (0,1)$ as the restriction $\RR\vert Y$ to a (any) subset $Y$ of measure $t$. The definition extends to any $t\in (0,\infty)$ by the formula $(\RR^{p})^{s}=\RR^{ps}$. This is easily shown to be consistent.

We now intend to extend such a definition to non ergodic equivalence relations. The reasonnable requirement is that $\RR^t$ ``meets equitably'' almost all the classes of $\RR$ or all the invariant Borel subsets. 
It will be defined for $t\in (0,\infty)$ when $t$ times the cardinal of the class $\vert \RR(x)\vert$ belongs to $\Nmath\cup\{\infty\}$, for almost every $x\in X$.

A Borel subset $B$ of $X$ is said to be \textbf{independent of $\RR$} if it is independent from the $\sigma$-algebra of $\RR$-invariant Borel subsets of $X$, i.e.  for any Borel subset $A\subset X$ that is a union of equivalence classes, one has $$\mu(A\cap B)=\mu(A)\mu(B).$$

\begin{proposition-def}[Contraction]\label{prop-def=amplification}
Let $\RR$ be a (non necessarily ergodic) p.m.p. countable standard equivalence relation on the standard atomless probability space $(X,\mu)$.
\begin{itemize}
\item[(a)] There exists a Borel subset $B$ that is independent of $\RR$ and has measure $\mu(B)=t$,\\
-- for every $t\in (0,1)$ when $\RR$ is aperiodic\\
-- for every rationnal number $t=\frac{p}{q}\in(0,1)$ $p,q\in \Nmath$ such that $q$ divides the cardinal of each class, otherwise\footnote{with the convention that any positive integer $q$ divides $\infty$}.

\item[(b)] If the borel subsets  $B$ and $C$ of $X$ are independent of $\RR$ and satisfy $\mu(B)=\mu(C)=t$, then there is a partial isomorphism $\varphi$ of the full groupoid $[[\RR]]$ such that $\varphi(B)=C$; in particular, the restrictions $\RR\vert B$ and $\RR\vert C$ are isomorphic. \\
This common isomorphism class is denoted by $\RR^{t}$ and called the $t$-\textbf{contraction} of $\RR$.\\
\item[(c)] More generally, for $s\in (0,1)\cup\Nmath^*$ and $t\in (0,1]$, the isomorphism class of $(\RR^s)^t$ depends only on the product $r=st$ as soon as it is defined. It is therefore denoted by $\RR^{r}$.
\end{itemize}
\end{proposition-def}

{\em Proof.}
(a) Take $(X,\mu)$ to be the interval $[0,1]$ with Lebesgue measure.
Consider the $\RR$-ergodic decomposition
$m:(X,\mu)\to \mathcal{EM}$ (where $\mathcal{EM}$ is the space of $\RR$-invariant ergodic measures) (see \cite{Var63} and consider $\RR$ as given by a -- non-necessarily free -- group action \cite[Th. 1]{FM77a}). Consider, for $t\in [0,1]$,  the Borel subset 
\begin{equation}X({t}):=\{x\in X: m(x)([0,x])\leq t\}.\end{equation}

When the classes are infinite, the invariant measures have no atom. It follows that $X({t})$ meets each ergodic component along a subset of measure $t$:
for each ergodic measure $e\in \mathcal{EM}$, $e(X({t}))= e(X({t})\cap m^{-1}(e))=t$. 

In case the classes are finite, an invariant ergodic measure is just the equiprobability on an equivalence class. 
The assertion follows easily.

(b) The proof is analoguous to the ergodic case. Let $G$ be a countable subgroup of the full group $[\RR]$ that generates $\RR$, i.e. for every $(x,y)\in \RR$, there is an element $g\in G$ such that $g(x)=y$ \cite[Th.~1]{FM77a}. Let $G=\{g_i\}_{i\in \Nmath}$ be an enumeration of $G$. Define $B_0:=B\cap g_0^{-1}(C)$, $C_0:=g_0(B_0)$ and recursively $B_{i+1}:=[B\setminus (B_0\cup B_1\cup\cdots\cup B_i)]\cap g_{i+1}^{-1}[B\setminus (C_0\cup C_1\cup\cdots\cup C_i)]$ and $C_{i+1}:=g_{i+1}(B_{i+1})$. Let $\varphi:\cup_{i\in \Nmath}B_i\to \cup_{i\in \Nmath}C_i$ be the partial isomorphism of $[[\RR]]$ defined by the restrictions $\varphi\vert B_i=g_i\vert B_i$. Consider the complement $B\setminus \cup_{i\in \Nmath}B_i$. By construction, its $\RR$-saturation $\bar{B}$ does not intersect $C\setminus \cup_{i\in \Nmath}C_i$. However, by independence, it intersects $B$ and $C$ equiprobably: $\mu(\bar{B}\cap B)=\mu(\bar{B}\cap C)=\mu(\bar{B}\cap \cup_{i\in \Nmath}C_i)=\mu(\varphi^{-1}(\bar{B}\cap \cup_{i\in \Nmath}C_i))=\mu(\bar{B}\cap \cup_{i\in \Nmath}B_i)$. From the equality of the extremal terms, it follows that $\bar{B}\cap B\setminus \cup_{i\in \Nmath}B_i$ is negligeable, and in turn that $\bar{B}$ itself is a null set.

Checking (c) is a routine calculation.
\hfill$\square$

We will use in several points the following useful observation. It renforce T{\"o}rnquist's theorem (Th.~\ref{th: S-i in general position}) by ``hiding any prescribed free action of a free group'' when putting two relations in general position.
For further applications, we don't want to require the relations $\SS_i$ to be ergodic.

\begin{lemma}[Replacement trick]\label{lem: replacement lemma-2} Let $\SS_1$, $\SS_2$ be two p.m.p. 
countable standard equivalence relations on the standard atomless probability space $(X,\mu)$.
Let $p=\frac{1}{t}\in \Nmath\setminus\{0,1\}$ be some integer such that both $\SS_i^{t}$ ($i=1,2$) make sense (Prop-Def. \ref{prop-def=amplification}).
\begin{itemize}
\item[(a)] If $\SS_1$, $\SS_2$ are in free product and there is some $Y\subset X$ with $\mu(Y)={t}$ for which both $\SS_i \vert Y\OE \SS_i^{t}$ ($i=1,2$), 
then 
\begin{equation}
(\SS_1*\SS_2)\vert Y= \SS_1 \vert Y *\SS_2\vert Y*\RR_{\alpha}
\end{equation}
where $\RR_{\alpha}$ is produced by some free action $\alpha$ of the free group $\FF_{p-1}$ on $Y$.

\item[(b)] Consider a free product decomposition
\begin{equation}\label{eq: replacement on Y}
\WW_1*\WW_2*\RR_{\beta} \textrm{\ \ \ with\ \ \ } \WW_i^{p}\OE \SS_i
\end{equation}
and $\RR_{\beta}$
is produced by SOME free action $\beta$ of the free group $\FF_{p-1}$ on $Y$.
Then it amplifies to a free product 
\begin{equation}\label{eq: amplification-of-the-free-product}
(\WW_1*\WW_2*\RR_{\beta})^{p}= \SS'_1 *\SS'_2 \textrm{\ \ \ where\ \ \ } \SS'_i\OE \SS_i.
\end{equation}

\item[(c)] For ANY free action $\beta$ of the free group $\FF_{p-1}$, there is a free product decomposition 
as in (\ref{eq: replacement on Y}). In particular, if $\beta$ has the relative property (T) then the resulting  $\SS'_1*\SS'_2$ in (\ref{eq: amplification-of-the-free-product}) also.
\end{itemize}
\end{lemma}

In case $\SS_1$ and $\SS_2$ are ergodic, any subset $Y$ of measure $\frac{1}{p}$ works and the assertion (a) is exactly Proposition~7.4~(2) in \cite{IPP05}.
The point (b) claims in particular that replacing $\alpha$ by any other free $\FF_{p-1}$-action leads to another general position between $\SS_1$ and $\SS_2$, hence the name.

{\em Proof of Lemma~\ref{lem: replacement lemma-2}}. (a) By assumption $\SS_{i}\simeq (\SS_{i}\vert Y)^p$. It follows that we have (for $i=1,2$) a free product decomposition $\SS_{i}=\SS_{i}\vert Y*\langle \varphi_{i,2}\rangle*\cdots*\langle \varphi_{i,p}\rangle$ where the graphings $(\varphi_{i,j}:Y\to Y_{i,j})_{j=2, \cdots, p}$ are such that each of the families $Y_{i,1}:=Y, Y_{i,2}, \cdots, Y_{i,p}$ ($i=1,2$) forms a partition of $X$. Thus,
\begin{equation}
\SS_{1}*\SS_{2}=\SS_{1}\vert Y*\SS_2\vert Y*\langle \varphi_{1,2}\rangle*\cdots*\langle \varphi_{1,p}\rangle*
\langle \varphi_{2,2}\rangle*\cdots*\langle \varphi_{2,p}\rangle
\end{equation}
and $(\varphi_{i,j})_{i=1,2; j=2,\cdots, p}$ is a treeing (see \cite{Gab00a}).
By considering all the compositions $\psi_{j,k}=(\varphi_{2,k})^{-1}\circ \varphi_{1,j}: \varphi_{1,j}^{-1}(Y_{1,j}\cap Y_{2,k})\to \varphi_{2,k}(Y_{1,j}\cap Y_{2,k})$ (for $j,k=2, \cdots, p$), 
we get another treeing $\Psi=(\psi_{j,k})_{j,k=2,\cdots,p}$ such that 
\begin{equation}
\SS_{1}*\SS_{2}=\SS_{1}\vert Y*\SS_2\vert Y*\langle \varphi_{1,2}\rangle*\cdots*\langle \varphi_{1,p}\rangle*
\langle \Psi\rangle
\end{equation}
and every point of $Y_{1,1}$ (resp. $Y_{2,1}$) belongs to the domains (resp. range) of $p-1$ partial isomorphisms of $\Psi$.
The measurable Hall's marriage theorem \cite[Th.~8]{Bollobas=1980} shows that there are $p-1$ isomorphisms $\tau_l:Y_{1,1}\to Y_{2,1}$ ($l=1, \cdots, p-1$) whose graphs (as subsets of $Y_{1,1}\times Y_{2,1}$) form a partitition of the union of the graphs of the $\psi_{j,k}$. In other words, up to subdivision, the $\psi_{j,k}$ can be assembled in $p-1$ isomorphisms $Y_{1,1}\to Y_{2,1}$  still forming a treeing.
Now remembering that $Y_{1,1}= Y_{2,1}=Y$, the map $a_l\mapsto \tau_l$ defines a free action $\alpha$ on $Y$ of the free group $\FF(a_1,\cdots, a_{p-1})$ so that 
\begin{eqnarray}
\SS_{1}*\SS_{2}&=&\SS_{1}\vert Y*\SS_2\vert Y*\langle \varphi_{1,2}\rangle*\cdots*\langle \varphi_{1,p}\rangle*
\langle \tau_1\rangle* \cdots *
\langle \tau_{p-1}\rangle\\
&=&\SS_{1}\vert Y*\SS_2\vert Y*\RR_{\alpha}*\langle \varphi_{1,2}\rangle*\cdots*\langle \varphi_{1,p}\rangle.
\end{eqnarray}
Since $(\varphi_{1,2}, \cdots,\varphi_{1,p})$ is a smooth treeing with fundamental domain $Y$, the assertion follows 
\begin{equation}
(\SS_{1}*\SS_{2})\vert Y=\SS_{1}\vert Y*\SS_2\vert Y*\RR_{\alpha}.
\end{equation}

(b) Let $a_2,\cdots, a_k, \cdots, a_{p}$ be a free generating set for the free group $\FF_{p-1}$ and $\beta(a_k):Y\to Y$ the associated isomorphisms. Let $\varphi_j:Y\times\{1\}\to Y\times \{j\}$ ($j=2, \cdots, p$) be the partial isomorphisms defined by the identity on the first coordinate. Thus on $Y\times \{1, 2, \cdots, p\}$, identifying $Y$ with $Y\times \{1\}$, (see eq.~(\ref{eq:amplification and free prod}))
\begin{eqnarray}
(\WW_1*\WW_2*\RR_{\beta})^{p}&=& (\WW_1*\WW_2*\RR_{\beta})*\langle \varphi_{2}\rangle*\cdots*\langle \varphi_{p}\rangle\\
&=& \WW_1*\WW_2*\langle \beta(a_2)\rangle*\cdots *\langle \beta(a_p)\rangle*\langle \varphi_{2}\rangle*\cdots*\langle \varphi_{p}\rangle\\
&=& \WW_1*\WW_2*\langle \varphi_{2}\circ\beta(a_2)\rangle*\cdots *\langle \varphi_{p}\circ\beta(a_p)\rangle*\langle \varphi_{2}\rangle*\cdots*\langle \varphi_{p}\rangle\\
&=& \underbrace{\WW_1*\langle \varphi_{2}\rangle*\cdots*\langle \varphi_{p}\rangle}_{\OE \WW_1^p\OE\SS_1}*\underbrace{\WW_2*\langle \varphi_{2}\circ\beta(a_2)\rangle*\cdots *\langle \varphi_{p}\circ\beta(a_p)\rangle}_{\OE \WW_2^p\OE\SS_2}
\end{eqnarray}

(c) Apply T{\"o}rnquist's theorem (Th.~\ref{th: S-i in general position}). 
It follows from 
\cite[Prop. 4.4 2$^{\circ}$, 4.5.3$^{\circ}$, 5.1, 5.2]{Pop06} that relative property (T) is inherited from a subrelation to a relation containing it, and is stable under amplification.
This concludes the proof of  Lemma~\ref{lem: replacement lemma-2}.
\hfill$\square$

\subsection{Proof of Th.~\ref{th: continuum many non vN rel. prescribed on factors}}
The general result follows from the case $\vert I\vert =2$ by splitting $(\SS_i)_{i\in I}$ into two families $(\SS_i)_{i\in I_1}$ and $(\SS_i)_{i\in I_2}$.
Take an integer $p\geq 4$. 
Consider now a free action $\beta$ of the free group $\FF_{p-1}$ such that one group element, say $a_2$, in some free generating set $(a_2,a_3, \cdots, a_{p})$ of $\FF_{p-1}$ acts ergodically, and the restriction to $\langle a_3, \cdots, a_{p}\rangle$ has relative property (T), for instance  
the standard free action on the $2$-torus $\Tmath^2=\Rmath^{2}/\Zmath^2$ of some $\FF_{p-1}< \SL(2,\Zmath)$ (see \cite{Bur91} and \cite[Lem.~6]{GP05}).
 In \cite[Prop. 1]{GP05} we constructed out of $\beta$ a continuum of free ergodic actions $\alpha_t$ of $\FF_{p-1}$, indexed by $r\in (0,1]$ defining a strictly increasing continuum of standard p.m.p. equivalence relations $\RR_{\alpha_{r}}$ such that for each $r\in (0,1]$, $\lim\limits_{s\to t, s<r}\nearrow\RR_{\alpha_{s}}=\RR_{\alpha_{r}}$, such that $\RR_{\alpha_{1}}=\RR_{\beta}$ and such that moreover the equivalence relation $\RR_0$ associated with the relative property (T) $\FF_{p-2}$-action $\alpha_0=\beta\vert \langle a_3, \cdots, a_{p}\rangle$ is contained in all the $\RR_{\alpha_{r}}$.
 
Apply Theorem~\ref{th: S-i in general position} to put $\SS_1^{\frac{1}{p}}, \SS_2^{\frac{1}{p}}$ and $\RR_{\beta}$ in general position as in Lemma~\ref{lem: replacement lemma-2}(c). 
\begin{equation}
\WW_1*\WW_2*\RR_{\beta} \textrm{\ \ \ with\ \ \ } \WW_i^{p}\OE \SS_i.
\end{equation}

 Since  $\RR_{\alpha_{r}}$ is a subrelation of $\RR_{\beta}$, it follows that $\WW_1$, $ \WW_2$, $ \RR_{\alpha_r}$ are freely independent and Lemma~\ref{lem: replacement lemma-2}(b) applies for every $r\in (0,1]$: $\WW_1*\WW_2*\RR_{\alpha_r}$ amplifies to 
\begin{eqnarray}
\SS_{r}:=(\WW_1*\WW_2*\RR_{\alpha_r})^p&=&\SS'_1*\SS'_2 \textrm{\ \ \ with\ \ \ } \SS'_i\OE\SS_i
\end{eqnarray}

 The family of factors $M_r=M(\SS_{r})$ associated with the equivalence relations $\SS_{r}$ thus satisfies the hypothesis of \cite[Prop. 4]{GP05}: it is a strictly increasing family of subfactors  $A\subset M_0\subset M_r\subset M_1$, with $A=L^{\infty}(X)$ Cartan, such that $M_1=\overline{\cup_{r<1} M_r}$ and $A\subset M_0$ relatively rigid. Then,  \\
(i) the pairs $A\subset M_t$ are relatively rigid Cartan subalgebras,\\
(ii) the classes of stable isomorphism among the $M_r$ are at most countable,
in particular there is an uncountable set $J$ such that for $r\in J\subset [0,1]$  the $M_r$  are pairwise not stably isomorphic.\\
Observe that moreover\\
(iii) $M_r$ has at most countable fundamental group (by \cite[Th. A.1]{NPS07})\\
(iv) $\SS_r$ has at most countable outer automorphism group (by \cite[Th. 4.4]{Pop06}).\\
The proof of th.~\ref{th: continuum many non vN rel. prescribed on factors}
is complete.\hfill$\square$

The above proof also show a version of Theorem~\ref{th: continuum many non vN rel. prescribed on factors}, with the same conclusion under the hypothesis that $\vert I\vert =2$, and the cardinal of the classes of $\SS_1$, $\SS_2$ admit a common divisor $p\geq 4$ (again with the convention that $p$ divides $\infty$).
In case $\SS_1$ is aperiodic and the classes of $\SS_2$ all have the same finite cardinality $p=2$ or $3$,
then put $\SS_1^{\frac{1}{n}}$, $\SS_2^{\frac{1}{p}}$ in general position by lemma~\ref{lem: replacement lemma-2}(c) for some $\FF_{p-1}$-action 
$\WW_1*\WW_2*\RR_{\beta} \textrm{\ \ \ with\ \ \ } \WW_i^{p}\OE \SS_i$
and then apply the above theorem to the free product decomposition $(\WW_1*\WW_2)*\RR_{\beta}$.
In case both $\SS_i$ are finite, $\SS'_1*\SS'_2$ can be made treeable and a specific treatment is applicable according to whether the result is amenable (cardinal of the classes are all $=2$) or not.

\bigskip
\paragraph{Proof of Theorem~\ref{th: continuum many non vN act prescribed on factors}}
If the groups are infinite, or if one can organize the $\Gamma_i$ in two families giving infinite groups $*_{i\in I_1} \Gamma_i$ and $*_{i\in I_2} \Gamma_i$, then Th.~\ref{th: continuum many non vN act prescribed on factors} follows from Th.~\ref{th: continuum many non vN rel. prescribed on factors}. Observe that if the $\SS'_i$ are orbit equivalent to relations that are produced by some free $\Gamma_i$-action, then the natural induced action of the free product $*_{i\in I} \Gamma_i$ is free when the $\SS'_i$ are in general position $*_{i\in I} \SS'_i$.
If the $\Gamma_i$ are all finite and $\vert I \vert =2$ or $3$, then $*_{i\in I} \Gamma_i$ is virtually a free group and the result follows from the free group study.
The remaining situation where  $\vert I\vert =2$, one of the $\Gamma_i$ is finite and the other is infinite follows from the above generalization of Theorem~\ref{th: continuum many non vN rel. prescribed on factors}.

\subsection{Proof of Theorem~\ref{th: cont many non vN rel. Out=1 prescribed on factors}}
It relies on the following:
\begin{theorem}\label{th: R-0 freprod R-Z Out=1}
Let $\RR_0$ be an ergodic relative property (T) p.m.p. equivalence relation on the non atomic standard probability space $(X,\mu)$.
There exist continuum many p.m.p. ergodic equivalence relations $(\RR_{t})_{t\in T}$, all contained in a standard equivalence relation $\overline{\RR}=\RR_0*\RR_{\sigma}$ for some free action $\sigma $ of $\FF_2$, such that
\begin{enumerate}
\item $\RR_{t}=\RR_0*\langle g_{t} \rangle$ for some isomorphism $g_{t}:X\to X$
\item the $\RR_{t}$ have relative property (T) 
\item the $\RR_{t}$ are pairwise von Neumann inequivalent 
\item $\Out(\RR_{t})=\{1\}$
\end{enumerate}
\end{theorem}
Assuming this (to be proved in the next section), we prove Theorem~\ref{th: cont many non vN rel. Out=1 prescribed on factors}: 
The general result follows from the case $\vert I\vert =2$ by splitting $(\SS_i)_{i\in I}$ into two families $(\SS_i)_{i\in I'}$ and $(\SS_i)_{i\in I''}$, and putting the $(\SS_i)_{i\in I'}$ (resp. $(\SS_i)_{i\in I''}$) mutually freely independent by theorem~\ref{th: S-i in general position}. 
Apply Lemma~\ref{lem: replacement lemma-2}(c) to $\SS_1,\SS_2$ with $p=4$, and a free action $\beta$ of $\FF_{3}=\ll a,b,c\rr$ whose restriction to $\FF_2=\ll a,b\rr$ has relative property (T) and is ergodic.
We get a free product decomposition:
$\WW_1*\WW_2*\RR_{\beta}$ with $\WW_i^{p}\OE \SS_i $. 
Further decompose $\RR_{\beta}=\RR_{\beta\ll a,b\rr}\freeprod \RR_{\beta\ll c\rr}$ and consider  the equivalence relation $\RR_0:=\WW_1\freeprod  \WW_2 \freeprod \RR_{\beta\ll a,b\rr}$ on $Y$.

Applying Theorem~\ref{th: R-0 freprod R-Z Out=1}, we get 
the family $\RR_{t}=\RR_0*\langle g_{t} \rangle= \WW_1\freeprod  \WW_2 \freeprod \RR_{\beta\ll a,b\rr}*\langle g_{t} \rangle$ for some isomorphism $g_{t}:Y\to Y$.
Now, $\RR_{\beta\ll a,b\rr}*\langle g_{t} \rangle$ may be interpreted as produced by a free action $\beta_{t}$ of $\FF_3$ to which Lemma~\ref{lem: replacement lemma-2}(b) applies.
\begin{eqnarray}
\SS_{t}:=\RR_t^p=(\WW_1\freeprod  \WW_2 \freeprod \RR_{\beta_t})^p
 &=&\SS'_{1,t}*\SS'_{2,t} \textrm{\ \ \ with\ \ \ } \SS'_{i,t}\OE\SS_i
\end{eqnarray}
The properties 2., 3. and 4. in theorem~\ref{th: R-0 freprod R-Z Out=1}, as well as ergodicity, being invariant under stable orbit equivalence of fixed amplification constant, it follows that the $\SS_t$ are ergodic, they have relative property (T), they are pairwise von Neumann inequivalent, and satisfy $\Out(\SS_{t})=\{1\}$. \hfill $\square$

\paragraph{Proof of Theorem~\ref{th: F3-action with trivial out}}
A stable isomorphism of relatively rigid von Neumann crossed-product $\FF_p\ltimes A$ implies a stable orbit equivalence \cite{Pop06} and the consideration of the first $\ell^2$-Betti number (see \cite{Gab02}) entails the triviality of the coupling constant.\hfill $\square$

\paragraph{Proof of Theorem~\ref{cor: free prod gps out=1}}
It follows immediately from Theorem~\ref{th: cont many non vN rel. Out=1 prescribed on factors} for infinite free factors. The same argument as in the proof of Theorem~\ref{th: continuum many non vN act prescribed on factors} reduces the case $\Gamma_1*\textrm{(finite group)}$ by using the
 generalized version of Theorem~\ref{th: continuum many non vN rel. prescribed on factors} (see the end of its proof) and the case $\textrm{(finite group)}*\textrm{(finite group)}*\textrm{(finite group)}$ to the study for the free group.
\hfill $\square$

\subsection{Proof of Theorem \ref{th: R-0 freprod R-Z Out=1}}

\begin{lemma}\label{lem: uncoutable implies equality somewhere}
Let $\RR_0, \RR$ be ergodic p.m.p. equivalence relations on the standard 
probability space $(X,\mu)$. Assume $\RR_0$ has relative property (T) and assume that $\RR_0\subset \RR$.
If $(\Delta_{t}:X\to X)_{t\in T}$ is an uncountable family of (measure preserving) automorphisms such that 
$\forall t\in T$,  $\Delta_{t}(\RR_0)\subset \RR$, then there exist $k\not= l\in T$ such that $\Delta_k=\Delta_{l}$ on a non-negligible Borel subset.
\end{lemma}
The $\Delta_{t}$ induce trace-preserving isomorphisms $\bigl(A\subset M(\RR_0)\bigr)\simeq \bigl(A\subset M(\Delta_{t}(\RR_0))\bigr)$ leading to embeddings $\theta_{t}:M(\RR_0)\subset M(\RR)$ with $\Delta_{t}(a)=a \circ  \Delta_{t}^{-1}$ for all $a\in A=L^{\infty}(X,\mu)$.
The Hilbert space 
$\HHH:=L^2(M(\RR),\tau)$ with its standard $M(\RR)-M(\RR)$-bimodule structure inherits for each $i,j\in T$ an $M(\RR_0)-M(\RR_0)$-bimodule structure $\HHH_{i,j}$, given by
$$ u \point\limits_{i} \xi \point\limits_{j} v= \Theta_i(u) \xi \Theta_{j}(v)$$
By separability of $\HHH$ and uncountability of $T$, the (tracial) vector $\xi$, image of $Id$ in the standard embedding $M(\RR)\subset L^{2}(M(\RR),\tau)$ satisfies for $\epsilon<1/2$ the condition of rigidity def.~\ref{def: rigidity} for $A\subset M(\RR_0)$,  for some $k\not=l$.
Thus there exists $\xi_0\in L^{2}(M(\RR), \tau)$ with 
\begin{equation}
\Vert \xi_0-\xi\Vert<\epsilon \label{eq: xi-0 near xi}
\end{equation}
 such that $\forall a\in A$, 
$a\point\limits_{k} \chi_0=\chi_0\point\limits_{l} a$, i.e. 
$\Theta_k(a) \chi_0=\chi_0\Theta_l(a)$.
By $A-A$-bimodularity of the projection $P:L^2(M(\RR),\tau)\to L^2(A)$, its image $h:=P(\xi_0)\in L^{2}(X,\mu)$ satisfies the same identities: $\forall a\in A=L^{\infty}(X)$, $ (a\circ \Delta_k^{-1}) h=(a\circ \Delta_l^{-1}) h$.
In particular, $\Delta_k^{-1}=\Delta_l^{-1}$ on the support of $h$ ($h\not=0$ by (\ref{eq: xi-0 near xi})).\hfill$\square$

\bigskip

We will now construct two treeings $\Phi=(\varphi_n)_{n\in \Nmath\setminus\{0\}}$ and $\Psi=(\psi_n)_{n\in \Nmath\setminus\{0\}}$. They will in particular satisfy: \\
-a- $\dom(\psi_n)=\dom(\varphi_n)$, $\rng(\psi_n)=\rng(\varphi_n)$, and both families $(\dom(\psi_n))_{n}$ and $(\rng(\psi_n))_n$ form a partition of $X$.\\
-b- They will be in general position, in the sense that : $\overline{\RR} :=  \ll \RR_0, \Phi, \Psi  \rr=\RR_0\freeprod \ll \Phi \rr \freeprod \ll \Psi \rr$.

Let's introduce a notation.
For each subset $E\subset \Nmath\setminus\{0\}$, let's denote $\Phi_{E}:=(\varphi_n)_{n\in E}$ and $\Psi_{\overline{E}}=(\psi_n)_{ n\not\in E}$.
The relations 
\begin{eqnarray}
\RRE &:=&\RR_0\freeprod \ll \Phi_E  \rr\\
\RRTE&:=&\RR_0\freeprod \ll \Phi_E \rr \freeprod \ll \Psi_{\overline{E}} \rr
\end{eqnarray}
will have relative property (T) since $\RR_0$ does.
By \cite[Th. 4.4]{Pop06} $\Aut(\RRTE)$ is countable modulo the full group $[\RRTE]$.
Observe that $ \ll \Phi_E \rr \freeprod \ll \Psi_{\overline{E}} \rr$ is a treed equivalence relation and that the partial isomorphisms $\Phi_E \cup \Psi_{\overline{E}}$ fit together to form a single automorphism $g_{E}:X\to X$.
\begin{equation}
\RRTE=\RR_0\freeprod \ll g_E \rr
\end{equation}
By Th.~\ref{th: S-i in general position}, choose a single automorphism $g_{\emptyset}$ of $X$ which is freely independent from $\RR_0$. Define the $\psi_n$ by restricting it to disjoint pieces of measures $\epsilon_n>0$ such that $\sum_{n\in \Nmath\setminus\{0\}}\epsilon_n=1$.

The $\varphi_n$'s are now constructed inductively:
for $n\in \Nmath$, let $\SS_n$ be the auxiliary p.m.p. standard equivalence relation generated by $\RR_0$, $\Phi_{\{1, 2, \cdots, n-1\}}$, $\Psi$ and all of the $\Aut(\RRTF)$ for $F\subset \{1, 2, \cdots, n\}$. The countability of $\Aut(\RRTF)$ modulo $[\RRTF]$ ensures that $\SS_n$ is actually countable. For $n=1$, simply define $\SS_1$ to be generated by $\RR_0$, $\Psi$ and $\Aut(\RR_0)$.
Choose an isomorphism $f_n$ of $X$ that is freely independent from $\SS_n$, restrict it to $\dom(\psi_n)$ and, by ergodicity of $\RR_0$, compose it by a certain $h\in [\RR_0]$ such that $h\circ f_n(\dom(\psi_n))=\rng(\psi_n)$. Then define $\varphi_n:=h\circ f_{n}{\vert \dom(\psi_n)}$.

\begin{theorem}\label{th: uncount. family R-E with trivial Out}
There is an uncountable family $\mathcal{E}$ of subsets of $\Nmath\setminus\{0\}$ such that for all $E\in \mathcal{E}$,  $\Out(\RRTE)$ is trivial for $E\in \mathcal{E}$ and the classes of von Neumann equivalence are at most countable.
\end{theorem}
The proof now follows the lines of \cite[Th. 4.1]{PV08}.
Choose an uncountable almost disjoint family $\mathcal{E}_1$ of $\Nmath\setminus\{0\}$ \cite{Sierpinski-1928}, i.e. a family of infinite subsets of $\Nmath\setminus\{0\}$ such that $E\cap F$ is finite for each $E\not= F$ in $\mathcal{E}_1$ and fix a $\Delta_E\in \Aut(\RRTE)$ for each $E\in \mathcal{E}_1$. We will show that at least one $\Delta_E$ is inner, i.e. belongs to $[\RRTE]$. It follows that the sub-family of those $\RRTE$ with non-trivial outer automorphisms group is
at most countable.

Step 1. There exists $E,F\in \EEE$ such that $E\not= F$ and $\Delta_E(x)=\Delta_F(x)$ for all $x$ in a non-negligible subset $U_{E,F}\subset X$.
We apply lemma~\ref{lem: uncoutable implies equality somewhere} to the family of $\Delta_E$ with $\Delta_E(\RR_0)\subset \RR_E\subset \overline{\RR}$.

Step 2. Since $\RR_0$ is ergodic and contained in $\RRTE\cap \RRTF$, there are $f_E\in [\RRTE]$ and $f_F\in [\RRTF]$ such that $\mu$-almost everywhere in $X$
\begin{equation}\label{eq: rel fe Delta fF}
f_E \Delta_E=f_F \Delta_F 
\end{equation}

Step 3. $\Delta_E$ belongs to $[\RRTE]$:
From the relation (eq. \ref{eq: rel fe Delta fF})
\begin{eqnarray}
\Delta:=f_E \Delta_E=f_F \Delta_F 
\ \ \in\ \  \Aut(\RRTE)\cap \Aut(\RRTF)\ \  \subset\ \  \Aut(\RRT_{E\cap F})
\end{eqnarray}

Let $n\in E$ be such that $n>\max(E\cap F)$. Since $\Delta\in \Aut(\RRTE)$ and $\varphi_n\in [[\RRTE]]$ the partial isomorphism 
\begin{equation}\label{eq: rel Delta-phi-n-h}
h:=\Delta \varphi_n \Delta^{-1}
\end{equation}
belongs to $[[\RRTE]]$.
Let $p$ be the smallest index for which $h\vert U\in [[\RR_0\freeprod \ll \Phi_{E\cap\{1, 2, \cdots, p\}}\rr \freeprod \ll \Psi_{\overline{E}} \rr]]$, for a non-negligeable domain $U$.

By the condition that the $\varphi_p$ are freely independent from $\SS_{q}$, for $q<p$, and since $\Delta\in \Aut(\RRT_{E\cap F})\subset [\SS_{\max(E\cap F)}]$, the relation (eq.~\ref{eq: rel Delta-phi-n-h}) 
entails that $p=n$.
Again by independence, by isolating the letters $\varphi_n^{\pm 1}$ in the reduced expression of $h$, the relation (eq.~\ref{eq: rel Delta-phi-n-h}) delivers a subword $w\in [[\RR_0\freeprod \ll \Phi_{E\cap\{1, 2, \cdots, p-1\}}\rr \freeprod \ll \Psi_{\overline{E}} \rr]]\subset [[\RRTE]]$ such that $\Delta=w$ on some non negligible set, i.e. $\Delta$ is ``locally inner''.
It follows by ergodicity of $\RR_0$ that $\Delta$ and thus that $\Delta_E$ belongs to $[\RRTE]$. 

\medskip
We show that the classes of orbit equivalence are at most countable.
For this, we show that
given $E\in \mathcal{E}_1$, there are at most countably many $F\in \mathcal{E}_1$ such that $\RRTE\OE\RRTF$. Apply lemma~\ref{lem: uncoutable implies equality somewhere} to an uncountable family of isomorphisms $\Delta_{F,E}: \RRTE\OE\RRTF$. Then at least two of them, for $F,F'\in\mathcal{E}_1$,  coincide on some non-negligible Borel subset. Like in step 2 above, it follows by ergodicity of $\RR_0\subset \RRTF\cap \RRT_{F'}$ that there are $f_F\in [\RRTF]$ and $f_{F'}\in [\RRT_{F'}]$ such that $\mu$-almost everywhere in $X$
\begin{equation}
  f_F\Delta_{F,E}=f_{F'}\Delta_{F',E} 
\end{equation}
It follows that $\RRTF=\RRT_{F'}$ and thus $F=F'$.

Now, von Neumann equivalence reduces to orbit equivalence.
The $A\subset M(\RRTE)$ are all rigid. An isomorphism $\Theta: M(\RRTE)\simeq M(\RRTF)$ 
delivers the relatively rigid Cartan subalgebra $\Theta(A)\subset M(\RRTF)=M(\RR_0)\freeprod\limits_{A} M(\ll g_E\rr)$. By \cite[7.12]{IPP05}, there is a unitary $u\in M(\RRTF)$ such that $u \Theta(A) u^* = A$, i.e. $u \Theta(.) u^*$ induces an orbit equivalence $\RRTE\OE \RRTF$. 
This completes the proof of theorem \ref{th: R-0 freprod R-Z Out=1}.\hfill $\square$

\bigskip

\noindent {\bf Acknowledgments} I am grateful to Sorin Popa and Stefaan Vaes for useful and entertaining discussions on these subjects of rigidity and von Neumann algebras. I would like to thank Julien Melleray and Mikaël Pichot too.

\def\cprime{$'$} \def\cprime{$'$} \def\cprime{$'$}

\nobreak
\noindent \textsc{D.~G.: CNRS - Université de Lyon,  ENS-Lyon, UMPA UMR 5669,
69364 Lyon
cedex 7, FRANCE}

\noindent \texttt{gaboriau@umpa.ens-lyon.fr}

\end{document}